\documentclass{amsart}

\usepackage{ifthen}
\usepackage[T1]{fontenc}
\usepackage{ae,aecompl}
\usepackage{amsmath,amsfonts,amssymb}
\usepackage{color}
\usepackage{dsfont} 

\usepackage[utf8]{inputenc}
\usepackage{graphicx}

\usepackage{setspace}
\usepackage[absolute,overlay]{textpos}
\usepackage{tikz}
\usepackage{pgfplots}
\usepgfplotslibrary{patchplots}
\usepgfplotslibrary{groupplots}
\pgfplotsset{compat=newest}
\usetikzlibrary{calc}
\usepackage{fancyhdr}
\usepackage{footnote}
\usepackage{todonotes}
\usepackage{nicefrac} 
\usepackage{bm}
\usepackage[boxruled]{algorithm2e}
\usepackage{comment} 

\newtheorem{theorem}{Theorem}

\newtheorem{remark}[theorem]{Remark}

\begin{document}

\title{Simulation-based High-Speed Elongational Rheometer for Carreau-type Materials}

\author[Kannengießer]{Lukas Kannengießer$^1$}

\author[Arne]{Walter Arne$^2$}

\author[Bier]{Alexander Bier$^3$}

\author[Marheineke]{Nicole Marheineke$^{1,\star}$}

\author[Schubert]{Dirk W. Schubert$^3$}

\author[Wegener]{Raimund Wegener$^2$}

\date{\today\\
$^1$ Universit\"at Trier, Arbeitsgruppe Modellierung und Numerik, Universit\"atsring 15, D-54296 Trier, Germany\\
$^2$ Fraunhofer ITWM, Fraunhofer Platz 1, D-67663 Kaiserslautern, Germany\\
$^3$ FAU Erlangen-Nürnberg, Lehrstuhl für Polymerwerkstoffe, Martensstr. 7, D-91058 Erlangen, Germany\\
$^\star$ corresponding author, email: marheineke@uni-trier.de, orcid: 0000-0002-5912-3465
}

\begin{abstract}
For the simulation-based design of fiber melt spinning processes, the accurate modeling of the processed polymer with regard to its material behavior is crucial. In this work, we develop a high-speed elongational rheometer for Carreau-type materials, making use of process simulations and fiber diameter measurements. The procedure is based on a unified formulation of the fiber spinning model for all material types (Newtonian and non-Newtonian), whose material laws are strictly monotone in the strain rate. The parametrically described material law for the elongational viscosity implies a nonlinear optimization problem for the parameter identification, for which we propose an efficient, robust gradient-based method. The work can be understood as a proof of concept, a generalization to other, more complex materials is possible.
\end{abstract}




\maketitle

\section{Introduction}
Melt spinning of polymers is an important process for the production of continuous fibers \cite{hufenus2020melt}, which is very often used in the technical textile industry, e.g., for the manufacture of hygiene products, filter media, insulation and soundproofing materials,  \cite{karian2003handbook}. Synthetic fibers currently dominate the industrial market due to their comparatively low production costs and wide range of applications. In the spinning process the molten polymer is forced through a fixed diameter die and then stretched either by a draw-off device or by aerodynamic stretching via an aspirator. In case of aerodynamic stretching, an air stream is directed into a tube so that it exerts a force on the fiber surface, which leads to drawing. Process design and fabric optimization ask for the ability to simulate the process. Thereby, a major challenge is the accurate modeling of the polymer fibers with respect to their material behavior. 

The elongational (extensional) viscosity relates the extensional stresses to the strain rates and was first investigated by Trouton in 1906 \cite{trouton1906coefficient}. He showed that in Newtonian fluids elongational and shear viscosity just differ by a factor of 3. However, polymers are not necessarily Newtonian fluids. In the hot spinning regime they might behave like generalized Newtonian fluids where the scalar-valued viscosity depends not only on temperature but also on strain rate. Various procedures have been developed to measure elongational viscosity, since it can be not deduced from shear viscosity information for generalized Newtonian and non-Newtonian fluids. The main disadvantage of all these devices is that they only work at elongational strain rates of up to $10\,\mathrm{s}^{-1}$. The strain rates occurring during fiber spinning can be at least one order of magnitude higher, so that the conventional measurement methods are not applicable  \cite{munstedt2018extensional}.

\begin{figure}[tb]
     \centering
     \includegraphics[scale=0.45]{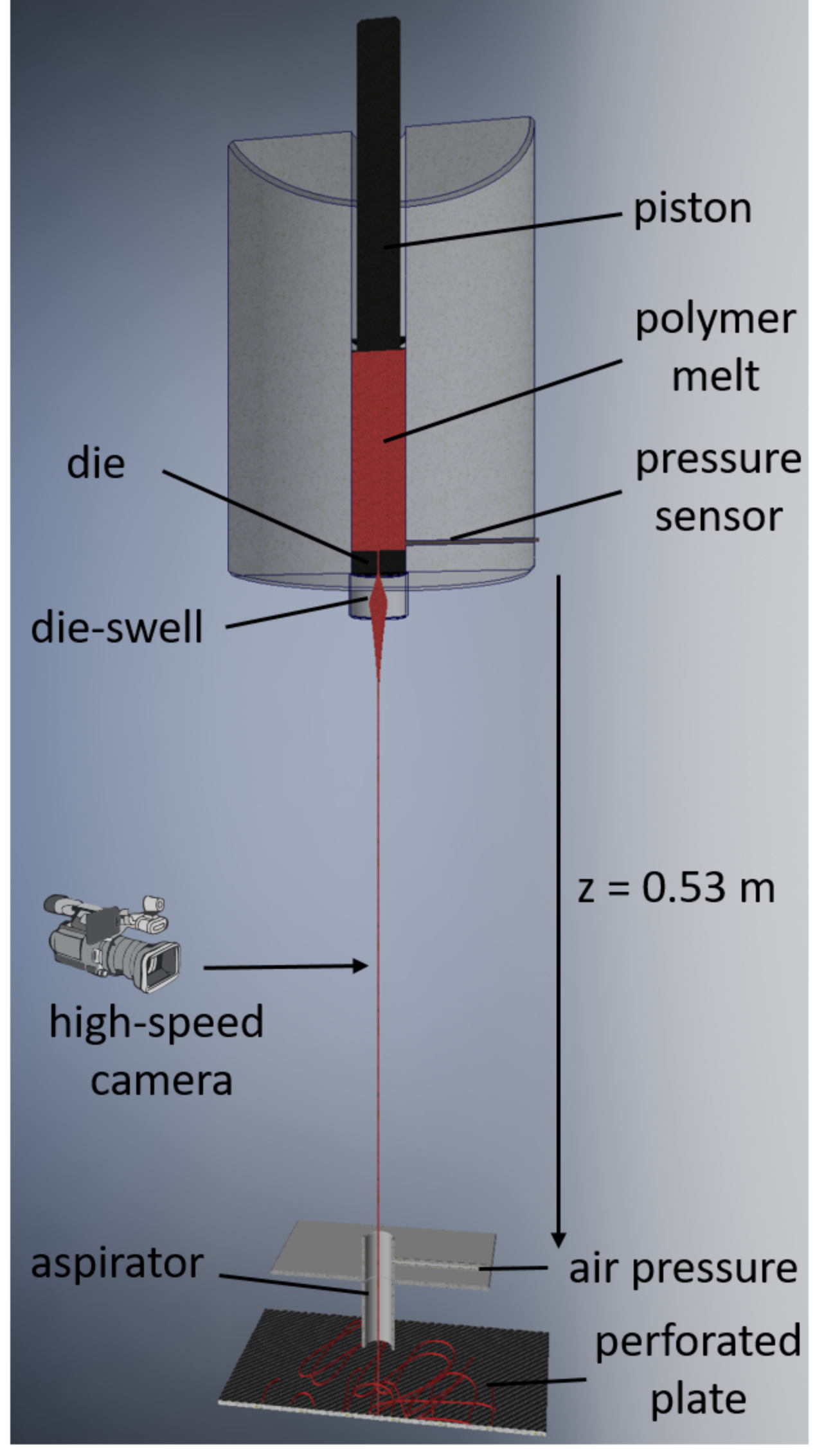}
     \caption{Spinning apparatus with single die and aspirator, used in \cite{bier2022novel}}
     \label{fig:prozessAbbildung}
 \end{figure}

Recently, a novel high-speed elongation rheometer has been introduced in \cite{bier2022novel}.  In a lab, a smaller spinning apparatus with a single die, correspondingly lower throughput and aspirator was considered, cf., Fig.~\ref{fig:prozessAbbildung}. The aspirator enables high speed and the generation of higher stretching rates. To avoid effects such as crystallization, fiber-air interactions an amorphous
polymethylmethacrylate (PMMA) with a linear structure was spun in still, isothermal air. On top of measurements for fiber diameter and shear viscosity, a rheometer was developed by help of aspirator theory \cite{kunzelmann2017korrelation}, \cite{qin2019simple}. Combining it with temperature approximations from Newtonian simulations, the findings suggested a non-Newtonian material behavior described by a Carreau-type model.
Concerning fiber spinning simulations, generalized Newtonian material models differ in the strain rate dependence of the elongational viscosity, whose handling in principle requires sophisticated material-dependent numerics. 

The aim of this work is the mathematical foundation of the high-speed elongational rheometer from \cite{bier2022novel} and its embedding in a closed, consistent simulation framework. 
We propose a procedure that allows for the unified formulation of the fiber spinning model as boundary value problem of explicit first order ordinary differential equations for all material types (Newtonian and non-Newtonian), whose material laws are strictly monotone in the strain rate. This provides a simplification and generalization of the numerical treatment and a reduction of the computational complexity.
For the class of generalized Newtonian fluids with Carreau-type material model we formulate the parameterization and determination of the elongational viscosity as a nonlinear optimization problem and present an efficient, robust gradient-based method for parameter identification. In contrast to \cite{bier2022novel}, we exclusively rely on the measurements and do without any further assumptions or heuristics. The study can be understood as a proof of concept for the proposed simulation-based parameter identification procedure. Application and generalization to other more complex materials is possible.

The article is structured as follows. Section~\ref{sec:model} deals with the fiber spinning simulation model. We introduce the unified formulation of the model equations for general material laws and present a collocation-continuation scheme for the numerical treatment. In Section~\ref{sec:param_ident} we explain the gradient-based parameter identification procedure. In Section~\ref{sec:results} we show and discuss the identified material parameters and fiber behavior for the spinning setup taken from \cite{bier2022novel} and investigate the performance of our approach. Details to model closure and simulations are provided in the appendix.

\section{Fiber Spinning Simulation Model}\label{sec:model}
In the following we present the fiber spinning model and discuss its numerical treatment. For the material behavior we consider a parameterically described material law of Carreau-type, which leads to a problem of parameter identification.
\subsection{Fiber model}\label{sec:FiberModel}
Uniaxial spinning in the direction of gravity can be described by a stationary one-dimensional fiber model based on the balances of mass, momentum and energy and on a material law with strain rate-dependent extensional viscosity. The state variables are the convective speed $u$, the contact force $N$ and the temperature $T$ of the fiber, as the mass balance implies a constant mass flow $Q$
\begin{align*}
\frac{\mathrm{d}}{\mathrm{d}s} (\rho A u) = 0  \qquad \rightarrow \qquad Q = \rho A u = \mathrm{const}
\end{align*} 
and thus, under the assumption of circular cross-sections, allows the elimination of a potential geometric state variable for cross-sectional area $A$ or diameter $d$, i.e.,
\begin{align}\label{eq:d}
A = A(u,\rho) = \frac{Q}{\rho u},\qquad d = d(u,\rho) = 2\sqrt{\frac{Q}{\pi\rho u}}.
\end{align}
The density $\rho$ is treated here as an argument to cover temperature dependencies. The model reads
\begin{align*}
Q \frac{\mathrm{d}}{\mathrm{d}s} u &= \frac{\mathrm{d}}{\mathrm{d}s} N + \rho A g + f_{air}\\
c_p Q \frac{\mathrm{d}}{\mathrm{d}s} T &= N \frac{\mathrm{d}}{\mathrm{d}s} u - \pi d \alpha (T-T_{air})\\
N &= A \mu_e (T, \frac{\mathrm{d}}{\mathrm{d}s} u) \frac{\mathrm{d}}{\mathrm{d}s} u
\end{align*}
with
\begin{align*}
\rho = \rho(T),\quad c_p = c_p(T),\quad f_{air} = f_{air}(u,d,s),\quad \alpha = \alpha(u,d,s).
\end{align*}
Fiber density $\rho$ and heat capacity $c_p$ are considered as temperature-dependent, the elongational viscosity $\mu_e$ as temperature- and strain rate-dependent. The air force per unit length $f_{air}$ and the heat transfer coefficient $\alpha$ depend not only on fiber velocity $u$ and diameter $d$, but also on the ambient airflow and  thus explicitly on the position $s$. The ambient temperature  $T_{air} = T_{air}(s)$ occurs also in the heat transfer. The gravitational acceleration is denoted by $g$.

We consider spinning processes on the interval $s \in [0, L]$, where $s = 0$ is identified with the nozzle position or a position near by (e.g., after a die swell not covered by the model) and $s=L$ with the take-off position or a take-up point. At the take-off we prescribe the fiber velocity or, alternatively, the diameter. This yields the boundary conditions
\begin{align*}
    T(0)=T_{in},\quad u(0)=u_{in},\quad \text{and} \quad u(L)=u_{out} \, \text{ or } \, d(u(L),\rho(T(L))=d_{out} .
\end{align*}

Each dimensional quantity $x$ is non-dimensionalized by help of a typical referential value $x_\circ$ according to $x=x_\circ \tilde{x}$ with dimensionless counterpart $\tilde{x}$. We leave open the exact choice of $s_\circ$, $u_\circ$, $T_\circ$, $\rho_\circ$, $c_{p,\circ}$, $\mu_\circ$, $\alpha_\circ$ and also $Q_\circ$, but for simplicity we choose
$L_\circ =s_\circ$, $N_\circ =Q_\circ u_\circ$, $f_\circ =Q_\circ u_\circ/s_\circ$, $u_{in,\circ} =u_{out,\circ} =u_\circ$, $T_{in,\circ} =T_\circ$. Since the auxiliary geometric quantities $A$ and $d$ can be inserted as state functions into the differential equations, it is useful to preserve their form in the non-dimensionalization, yielding
$A_\circ = Q_\circ/(\rho_\circ u_\circ)$ and $d_\circ=\sqrt{A_\circ}$. The only remaining model constant $g$ merges into the dimensionless Froude number.
Dropping the superscript $\tilde{\,}$ for readability, the dimensionless system is given by, $s\in[0,L]$, 
\begin{subequations}\label{eq:sys}
\begin{align} \label{sys:1a}
Q\frac{\mathrm{d}}{\mathrm{d}s} u &= \frac{\mathrm{d}}{\mathrm{d}s} N + \frac{Q}{\mathrm{Fr}^2}\,\frac{1}{u} + f_{air}\\ \label{sys:1b}
c_p Q \frac{\mathrm{d}}{\mathrm{d}s} T &= \mathrm{Ec} \, N \frac{\mathrm{d}}{\mathrm{d}s} u - \mathrm{St}\, \pi d \alpha (T-T_{air})\\ \label{sys:1c}
\frac{\mathrm{Re}}{Q}\,\rho\, u N &= \mu_e (T, \frac{\mathrm{d}}{\mathrm{d}s} u) \frac{\mathrm{d}}{\mathrm{d}s} u
\end{align}
with 
\begin{align*}
\rho = \rho(T),\quad c_p = c_p(T),\quad f_{air} = f_{air}(u,d,s),\quad \alpha = \alpha(u,d,s), \quad
d = d(u,\rho) = 2\sqrt{\frac{Q}{\pi\rho u}},
\end{align*}
supplemented by
\begin{align}
 T(0)=T_{in},\quad u(0)=u_{in},\quad \text{and} \quad u(L)=u_{out} \, \text{ or } \, d(u(L),\rho(T(L))=d_{out} 
\end{align}
\end{subequations}
The dimensionless numbers are the Reynolds number $\mathrm{Re}$ (ratio of inertial and viscous forces), the Froude number $\mathrm{Fr}$ (ratio of inertial and gravitational forces), the Eckert number $\mathrm{Ec}$ (ratio of the kinetic energy and enthalpy difference) and the Stanton number $\mathrm{St}$ (measure of the relative cooling intensity during heat transfer),
\begin{align*}
\mathrm{Re} = \frac{\rho_\circ u_\circ s_\circ}{\mu_\circ},\qquad \mathrm{Fr} = \frac{u_\circ}{\sqrt{g L_\circ}},\qquad \mathrm{Ec} = \frac{u_\circ^2}{c_{p,\circ}T_\circ},\qquad \mathrm{St} = \frac{\alpha_\circ d_\circ s_\circ}{c_{p,\circ}Q_\circ}
\end{align*}

\begin{remark}
Note that the model equations can be simplified by
clever choice of typical quantities as model constants, e.g.,
\begin{itemize}
    \item choosing $Q_\circ$ as the constant mass flux of the problem yields the dimensionless $Q=1$ 
    \item choosing $L_\circ = s_\circ$ as the spinning length yields the dimensionless $L=1$ and thus $s\in[0,1]$
    \item choosing $T_\circ$ and $u_\circ$ as the temperature and velocity at the nozzle, the boundary conditions simplify to $T(0) = 1$ and $u(0) = 1$ as well as $u(L) = \mathrm{Dr}$ with the dimensionless draw ratio $\mathrm{Dr}$  (ratio of take-off and spin velocity)
\end{itemize}
The simplifications make sense when dealing with a single problem. But if one aims at a common non-dimensionalization for various scenarios, which differ in the model constants, as it is the case in optimization, the more general approach should be used.  
\end{remark}

We assume that the material law $P_e(T,\frac{\mathrm{d}}{\mathrm{d}s} u)=\mu_e (T, \frac{\mathrm{d}}{\mathrm{d}s} u) \frac{\mathrm{d}}{\mathrm{d}s} u$ is strictly monotone in the strain rate. Then, its inverse $P^{-1}_e(T,.)$ exists and the fiber model can be formulated as boundary value problem of explicit first order ordinary differential equations. The constitutive equation \eqref{sys:1c} particularly becomes the differential equation for the velocity $u$,
$$\frac{\mathrm{d}}{\mathrm{d}s} u=P^{-1}_e(T,\frac{\mathrm{Re}}{Q}\rho(T)uN).$$
If the inverse is not analytically available, we introduce the strain rate as additional variable $\dot{\epsilon}$, differentiate the constitutive equation \eqref{sys:1c} and use the resulting equation as differential equation for the strain rate. As boundary condition we pose the constitutive equation at the nozzle, 
\begin{subequations}\label{eq:dot-eps}
\begin{align}
    \frac{\mathrm{d}}{\mathrm{d}s} u&=\dot{\epsilon}\\
    \frac{\mathrm{d}}{\mathrm{d}s} \dot{\epsilon}&=
    \frac{ \frac{\mathrm{Re}}{Q}(\rho N \dot{\epsilon}+\rho u \frac{\mathrm{d}}{\mathrm{d}s} N)+(\frac{\mathrm{Re}}{Q}uN \partial_T \rho -\dot{\epsilon}\partial_T \mu_e)\frac{\mathrm{d}}{\mathrm{d}s}T}{\mu_e+\partial_{\dot{\epsilon}}\mu_e\dot{\epsilon}}, \qquad
    \frac{\mathrm{Re}}{Q}\rho(T)uN-\mu_e(T,\dot{\epsilon})\dot{\epsilon}\bigg|_{s=0}=0.
\end{align}
\end{subequations}
Note that the positivity of the denominator is hereby ensured by the strict monotony of $P_e(T,.)$. 

The broad class of generalized Newtonian fluids satisfies the monotonicity assumption and allows the (re-)formulation \eqref{eq:dot-eps}, cf., \cite{kannengiesserECMI}.
The Carreau fluid is a type of generalized Newtonian fluid. Its elongational viscosity in dimensionless form is modeled as
\begin{align*} 
\mu_e(T,\dot{\epsilon})=\mu_{e,\infty}(T)+ (\mu_{e,0}(T)-\mu_{e,\infty}(T))\left(1+\mathrm{De}^2\left( \lambda(T)\dot{\epsilon}\right)^2\right)^{\frac{n-1}{2}}
\end{align*}
with power index $n\geq0$, temperature-dependent relaxation time $\lambda$ as well as zero and infinite strain-rate viscosities $\mu_{e,0}$, $\mu_{e,\infty}$ and Deborah number $\mathrm{De}$. In the strain thinning regime, $n< 1$, of interest, the Carreau model can be considered as interpolation between Newtonian models with $\mu_{e,0}$ at $\dot\epsilon= 0$ and $\mu_{e,\infty}$ for $\dot\epsilon\rightarrow \infty$, in particular it holds $\mu_{e,\infty}(\cdot)=\lim_{\dot\epsilon \rightarrow \infty} \mu_{e}(\cdot,\dot\epsilon) \ll \mu_{e,0}(\cdot)$. Expecting moderate strain rates in the application, a simplified Carreau-type model has been proposed in \cite{bier2022novel}, setting $\mu_{e,\infty}(T)=0$ and $\lambda(T)=\mu_{e,0}(T)/K$. The simplified model hence depends only on two scalar-valued parameters $n$, $K$ and the temperature-dependent zero strain-rate elongational viscosity $\mu_{e,0}$.
For $n \rightarrow 1$ or $K \rightarrow \infty$ it provides Newtonian behavior, $\mu_e(\cdot, \cdot) \rightarrow \mu_{e,0}(\cdot)$. Elongational and shear viscosity are proportional, i.e., $\mu_e=\mathrm{Tr}\,\mu_s$. In generalized Newtonian materials, the Trouton number is $\mathrm{Tr}=3$. Thus, $\mu_{e,0}$ can be prescribed in terms of the zero-shear rate viscosity $\mu_{s,0}$, e.g., by help of a Vogel-Fulcher-Tammann model with material constants $\mu_c$, $B$ and Vogel temperature $T_{VF}$. In this work we follow \cite{bier2022novel} and consider the five-parametric model for the elongational viscosity,
\begin{subequations}\label{eq:CarreauLikeLawK}
\begin{align} \label{eq:CarreauLikeLaw}
\mu_e(T,\dot{\epsilon})&=\mu_{e,0}(T)\left(1+\mathrm{De}^2\left( \frac{\mu_{e,0}(T)}{K}\dot{\epsilon}\right)^2\right)^{\frac{n-1}{2}}\\\nonumber
    \mu_{e,0}(T)&=3\mu_{s,0}(T), \qquad \mu_{s,0}(T)=\mu_c \exp\left(\dfrac{B}{T-T_{VF}}\right).
\end{align}
Here, $\mu_c$, $B$, $T_{VF}$ can be straightforward determined from shear viscosity measurements for small shear rates. For mathematical reasons (cf., Sec.~\ref{subsec:opt-alg}) we introduce $\kappa$ via
\begin{align}
    K=\exp(\kappa).
\end{align}
\end{subequations}
Then, the identification of the following two scalar-valued material parameters remains
\begin{align}
\mathbf{p}=(n,\kappa)\in \Omega_p\subseteq[0,1]\times \mathbb{R}
\end{align}
Note that the inverse $P_e^{-1}(T,\cdot)$ of the resulting material law is not analytically given.

\begin{remark}\label{rem:K}
In general, the dimensionless Deborah number characterizes the rheology by representing the ratio of the relaxation and observation times, i.e., $\mathrm{De}=\lambda_\circ u_\circ/ s_\circ$. In the simplified Carreau-type model \eqref{eq:CarreauLikeLaw}, $\lambda_\circ$ can be traced to $\mu_\circ$ and $K_\circ$, i.e., $\lambda_\circ=\mu_\circ/K_\circ$. Thus, since $K$ is subject of our parameter identification, the Deborah number can here be viewed as a simple scaling factor. It can be arbitrarily set due to the freedom in $K_\circ$. 
\end{remark}

\subsection{Collocation-continuation scheme}\label{subsec:bvp_colloc}
Formulating the fiber spinning model as boundary problem of explicit first-order differential equations, we use a collocation scheme for the numerical approximation and solve the resulting nonlinear system with Newton method. Thereby, we compute the Jacobian by means of differentiation with complex variables. To provide suitable initializations for the Newton method and ensure global convergence we embed the collocation into a continuation / homotopy method. In the homotopy method we consider a family of boundary value problems with continuation parameter $c\in[0,1]$. Proceeding from the solution of a simple auxiliary problem for $c=0$, we follow a continuation path and solve a sequence of problems to finally obtain the solution for our original fiber spinning problem at $c=1$.

The family of boundary value problems for the state $\boldsymbol{y}=(u,N,T,\dot \epsilon)^T$ on $s\in[0,L]$ with material model parameter $\mathbf{p}=(n,\kappa)$ and continuation parameter $c$ reads
\begin{align}\label{eq:num_c}
    \frac{\mathrm{d}}{\mathrm{d}s} \boldsymbol{y}&=\boldsymbol{f}(\boldsymbol{y};\mathbf{p};c), \qquad
    \boldsymbol{g}(\boldsymbol{y}(0),\boldsymbol{y}(L);\mathbf{p};c)=\mathbf{0}
\end{align}
where
\begin{align*}
&\boldsymbol{f}(\boldsymbol{y};\mathbf{p};c) = \begin{pmatrix}
             \dot{\epsilon}\\
            f_N\\
         f_T    \\  
         \dfrac{ \frac{(c \mathrm{Re}+(1-c) \mathrm{Re}_0)}  {Q}(\rho N \dot{\epsilon}+\rho u f_N+uN \partial_T \rho f_T)
         -\dot{\epsilon}\partial_T \mu_e(\mathbf{p}) f_T}{\mu_e(\mathbf{p})+\dot{\epsilon}\partial_{\dot{\epsilon}}\mu_e(\mathbf{p})}
         \end{pmatrix}\\[0.2em]
&\boldsymbol{g}(\boldsymbol{y}(0),\boldsymbol{y}(L);\mathbf{p};c)=
\begin{pmatrix}
        u(0) - u_{in}\\
        d(u(L),\rho(T(L)) - (c d_{out}+(1-c) d(u_{in},\rho(T_{in}))\\
        T(0) - T_{in} \\ 
        \frac{(c \mathrm{Re}+(1-c) \mathrm{Re}_0)}  {Q}\rho(T_{in})u_{in} N(0)-\dot{\epsilon}(0)\mu_e(T_{in},\dot{\epsilon}(0);\mathbf{p})
               \end{pmatrix}
    \end{align*}
with abbreviations 
\begin{align*}
    f_N=Q\dot{\epsilon} - c\left( f_{air} +\frac{Q}{\mathrm{Fr}^2u} \right), \qquad f_T=\frac{c}{c_pQ}(\mathrm{Ec} N \dot{\epsilon}-\mathrm{St} \pi d\alpha(T-T_{air})). 
\end{align*} 
The isothermal stress-free fiber, i.e., $u\equiv u_{in}$, $N\equiv 0$, $T\equiv T_{in}$, $\dot\epsilon\equiv 0$, solves the auxiliary problem with $c=0$ for every $\mathbf{p}\in \Omega_p$ and acts as initialization for the continuation. Here, $\mathrm{Re}_0$ serves as initial Reynolds number. Our fiber spinning problem corresponds to $c=1$. Details to the step size strategy of the continuation used in the simulation can be found in Appendix~\ref{sec: NumericalBVPs}.
For the discretization of the boundary value problem we use a three-stage Lobatto IIIa formula. The collocation polynomial provides a once continuously differentiable solution that is fourth-order accurate uniformly in the interval of integration. Mesh selection and error control are based on the residual of the continuous solution, \cite{kierzenka2001bvp}.

\subsection{Parameter dependence}
In the following we assume that our fiber spinning model $\boldsymbol{\mathcal{S}}$ is well-posed for every material parameter $\mathbf{p}\in \Omega_p \subset \mathbb{R}^2$, $\Omega_p$ compact. Then, there exists a unique mapping between material parameter and fiber state
\begin{align}\label{eq:y-p-map}
    \boldsymbol{y}:\Omega_p \rightarrow (\mathcal{C}^1([0,L]))^4 \qquad \text{ with } \quad \boldsymbol{\mathcal{S}}(\boldsymbol{y}(\mathbf{p}),\mathbf{p})=\mathbf{0}
\end{align}
given by \eqref{eq:sys}-\eqref{eq:CarreauLikeLawK}.
On the discrete level we have the analogon
\begin{align}\label{eq:y-p-map-disc}
    \mathbf{y}:\Omega_p \rightarrow \mathbb{R}^{4n_c} \qquad \text{ with } \quad \boldsymbol{{S}}(\mathbf{y}(\mathbf{p}),\mathbf{p})=\mathbf{0}.
\end{align}
The nonlinear system $\boldsymbol{{S}}(
\mathbf{y}(\mathbf{p}),\mathbf{p})=\mathbf{0}$ consists of the parameter-dependent collocation equations for the state values on the grid $\triangle_c$, i.e., $\mathbf{y}(\mathbf{p})\in \mathbb{R}^{4n_c}$ with $|\triangle_c|=n_c$. The function $\boldsymbol{S}:\mathbb{R}^{4n_c} \times \Omega_p  \rightarrow \mathbb{R}^{4n_c}$ is differentiable in $\mathbf{y}$ and $\mathbf{p}$.
To solve the nonlinear system for fixed $\mathbf{p}$, we apply the Newton method for which we determine the Jacobian $\partial_y \boldsymbol{S}$ by means of differentiation with complex variables, \cite{squire1998}. Note that the regularity of the Jacobian yields the local differentiability of $\mathbf{y}$ with respect to $\mathbf{p}$ according to the Implicit Function Theorem, i.e., \begin{align}\label{eq:y-deriv}
\mathrm{D}_p \mathbf{y}(\mathbf{p})=(\partial_y \boldsymbol{S}(\mathbf{y}(\mathbf{p}),\mathbf{p}))^{-1} \,\partial_p \boldsymbol{S}(\mathbf{y}(\mathbf{p}),\mathbf{p}) .
\end{align}

\section{Gradient-based Parameter Identification}\label{sec:param_ident}

The approach for a novel high-speed elongation rheometer by \cite{bier2022novel} provides estimates for the material parameters based on measurements of the fiber diameter and further assumptions and heuristics. Our parameter identification makes only use of the measured data.

\subsection{Problem formulation}\label{subsec: opt}
Let $\mathbf{d}^{(k)}_\text{meas}\in \mathbb{R}^{n_m^{(k)}}$ be the data vector of measured diameters in the $k$th experiment, $k=1,\dots, M$, that is made dimensionless just like the respective model quantity.
Our idea is the identification of the material parameter $\mathbf{p}\in \Omega_p$ for the Carreau-like viscosity model \eqref{eq:CarreauLikeLawK} by solving a weighted least-squares problem of the form,
\begin{align*}
\min_{\mathbf{p}\in \Omega_p}\, \sum_{k=1}^M \| \mathbf{d}(\boldsymbol{y}^{(k)}(\mathbf{p})) - \mathbf{d}^{(k)}_\text{meas}\|^2_{\mathbf{W}(\boldsymbol{y}^{(k)}(\mathbf{p}))}.
\end{align*}
Here, $$\mathbf{d}(\boldsymbol{y}^{(k)})\in \mathbb{R}^{n_m^{(k)}},\qquad d_i(\boldsymbol{y}^{(k)})=d(u^{(k)}(s_i),\rho(T^{(k)}(s_i))), \quad i=1,\dots,{n_m^{(k)}}$$
are the diameters obtained from our fiber spinning model \eqref{eq:y-p-map} at the measurement points of the $k$th experiment, cf., \eqref{eq:d}. The considered norm relies on a symmetric positive definite fiber state-dependent weighting matrix $\mathbf{W}(\boldsymbol{y}^{(k)})\in \mathbb{R}^{{n_m^{(k)}}\times {n_m^{(k)}}}$, i.e., $\|\mathbf{x}\|^2_\mathbf{W}=\langle \mathbf{x},\mathbf{W}\mathbf{x}\rangle$ with Euclidean scalar product $\langle. \,,.\rangle$. In the following we discuss the data and specify the weighting strategy.

\subsection{Data and weighting strategies}
In the spinning process the fiber diameter obeys a monotonically decreasing behavior along the spinline. However, the measurement data are subject to inaccuracies. In particular, measurements near the take-up point are susceptible to interference.
Hence, we smooth the data with respect to an ansatz function.

We make use of the fact that the inverse squared diameter behaves like a velocity and consider the following velocity ansatz function proposed in \cite{bier2022novel},
\begin{align*} 
            u_f(s;b,c,v,u_0) = \dfrac{u_0\, v \exp \left(\left(\frac{s}{c} \right)^b\right)}
            {(\exp \left(\left(\frac{s}{c} \right)^b\right)-1)u_0+v}
\end{align*}
Since $u_f(0;\,.\,)=u_0$ holds, for each measurement series (spun fiber) $k\in \{1,\dots, M\}$ we set $u_0=(d^{(k)}_{\text{meas}, in})^{-2}$ (inverse squared diameter at or near the nozzle) and determine the remaining parameters $(b,c,v)\in \mathbb{R}^3$ by means of a nonlinear least-squares fit with respect to the (converted) data $\{(d^{(k)}_{\text{meas},s_i})^{-2},\, i=1,...,{n_m^{(k)}}\}$ and the relative deviation. Then, we use the smooth diameter profile
\begin{align}\label{eq:ansatzFunc}
d^{(k)}_\text{fit}(s)&=\sqrt{\dfrac{1}{u_f(s;({b},{c},{v})^{(k)},u_0^{(k)})}},\\ \nonumber
&\quad ({b},{c},{v})^{(k)}=\text{arg}\min_{(b, c, v) \in \mathbb R^3} \sum_{i=1}^{{n_m^{(k)}}}[u_f(s_i;b,c,v,u_0^{(k)})(d^{(k)}_{\text{meas},s_i})^{2} - 1]^2
\end{align}
for the identification of the material parameters.

Non-Newtonian material behavior becomes evident at moderate/large strain rates, and, in this sense, the Carreau-like model can be interpreted as a modification of the Newtonian behavior for moderate/large strain rates.
In the parameter identification, we emphasize regions with pronounced strain rates and stronger non-Newtonian behavior and hence weight the data with respect to the strain rates.
The optimization problem on the continuous level becomes
\begin{align}
\min_{\mathbf{p}\in \Omega_p} \mathcal{J}(\mathbf{p}),
\qquad \mathcal{J}(\mathbf{p})= \sum_{k=1}^M \| \sqrt{\dot{\epsilon}^{(k)}(\mathbf{p})}\,\,[d(u^{(k)}(\mathbf{p}),\rho(T^{(k)}(\mathbf{p})) - {d}^{(k)}_\text{fit}]\|^2_{\mathcal{L}^2([0,L^{(k)}])}
\end{align}
with fiber state $\boldsymbol{y}^{(k)}=(u,N,T,\dot\epsilon)^{(k)}$ of \eqref{eq:y-p-map} and the Lebesque norm $\|x\|_{\mathcal{L}^2([0,1])}^2=\int_0^1 x(s)^2\, \mathrm{d}s$. As discrete counter part we have
\begin{align}\label{eq:opt:disc}
\min_{\mathbf{p}\in \Omega_p} J(\mathbf{p}),
\qquad J(\mathbf{p})= \sum_{k=1}^M \| \mathrm{diag}(\sqrt{\boldsymbol{\Delta}^{(k)}})\,\mathrm{diag}(\sqrt{\boldsymbol{\dot{\epsilon}}^{(k)}(\mathbf{p})})\,\,[\mathbf{d}(\mathbf{y}^{(k)}(\mathbf{p})) - {\mathbf{d}}^{(k)}_\text{fit}]\|^2_{2}
\end{align}
in the Euclidian norm $\|\cdot\|_2$ with $\mathbf{y}^{(k)}\in \mathbb{R}^{4n_{c}^{(k)}}$ solution vector of the collocation \eqref{eq:y-p-map-disc} and
$\mathbf{d}(\mathbf{y}^{(k)})$, $\boldsymbol{\dot{\epsilon}}^{(k)}$, $\mathbf{d}_\text{fit}^{(k)}  \in \mathbb{R}^{n_{o}^{(k)}}$ vectors of respective quantities on the optimization grid $\triangle_{o}^{(k)}$, $|\triangle_{o}^{(k)}|=n_{o}^{(k)}$. We choose the optimization grid as subset of the collocation grid to avoid additional interpolation. The vector $\boldsymbol{\Delta}^{(k)}$ contains grid size information for the quadrature, $\Delta_i^{(k)}=s_i-s_{i-1}$, $i=1,...,n_o^{(k)}$. The cost functional $J$ can be also expressed in terms of a weighted Euclidian norm on the diameter information (cf., Sec.~\ref{subsec: opt}), then the symmetric positive definite weighting matrix is the diagonal matrix $\mathbf{W}(\mathbf{y}^{(k)})=\mathrm{diag}(\boldsymbol{\Delta}^{(k)})\,\mathrm{diag}(\boldsymbol{\dot{\epsilon}}^{(k)})$.

In the following we aim at solving \eqref{eq:opt:disc} being a nonlinear least squares problem of the form
\begin{align}\label{eq:opt}
\min_{\mathbf{p}\in \Omega_p} J(\mathbf{p}), \qquad J(\mathbf{p})=\|\mathbf{F}(\mathbf{Y}(\mathbf{p}))\|_2^2, \qquad \mathbf{Y}=(\mathbf{y}^{(1)},\ldots,\mathbf{y}^{(M)}):\Omega_p\rightarrow \mathbb{R}^{4\sum_k n_c^{(k)}}
\end{align}
with continuously differentiable function $\mathbf{F}:\mathbb{R}^{4\sum_k n_c^{(k)}}\rightarrow\mathbb{R}^{\sum_k n_o^{(k)}}$.
The differentiability of the cost function $J$ with respect to the material parameter $\mathbf{p}$ follows directly from the differentiability of $\mathbf{Y}$, cf., \eqref{eq:y-deriv}, and since $\Omega_p$ is compact the minimum exists.

\subsection{Algorithmic procedure}\label{subsec:opt-alg}
A function evaluation of $J$ in \eqref{eq:opt} is computationally expensive, since the solving of $M$ boundary value problems is necessary to determine the discretized fiber states $\{\mathbf{y}^{(k)}(\mathbf{p}),\, k=1,...,M\}$. Solving a fiber spinning problem is in general sophisticated, as the performance of the underlying Newton method crucially depends on a good initialization. That is why we use a continuation framework for the simulation, cf., Sec.~\ref{subsec:bvp_colloc}. However, in the context of parameter identification, we can dispense the continuation routine and instead make use of the iteration procedure in the optimization.
Having computed the fiber state $\mathbf{Y}(\mathbf{p}_l)$ at material parameter $\mathbf{p}_l$, we can use this result as initial guess for computing $\mathbf{Y}(\mathbf{p}_{l+1})$, if $ \Vert\mathbf{p}_{l+1} -\mathbf {p}_l \Vert_2$ is small enough.
Thus, we apply a trust-region Gauss-Newton method equipped with a suitable upper bound for the allowed trust region in order to not exceed this range, \cite{coleman1996interior}.

The optimization procedure relies on information for first and second derivatives. The gradient $\nabla J = 2 (\mathrm{D}_y\mathbf{F}\,\mathrm{D}_p\mathbf{Y})^T \,\mathbf{F} $ and the used approximation of the Hessian $\nabla^2 J  \approx 2 (\mathrm{D}_y\mathbf{F}\,\mathrm{D}_p\mathbf{Y})^T (\mathrm{D}_y\mathbf{F}\,\mathrm{D}_p\mathbf{Y})$ depend on the Jacobian $\mathrm{D}_y\mathbf{F}$ of $\mathbf{F}$ and on the Jacobian $\mathrm{D}_p\mathbf{Y}$ of $\mathbf{Y}$.
Whereas $\mathrm{D}_y\mathbf{F}$ is analytically available, this is not the case for $\mathrm{D}_p\mathbf{Y}$. An approximation via finite differences is not suitable because of the expensive function evaluations of $\mathbf{Y}$. Instead, we use the relation from \eqref{eq:y-deriv} for the numerical computation; i.e., $\boldsymbol{S}^{(k)}(\mathbf{y}^{(k)}(\mathbf {p}),\mathbf{p}) = \mathbf{0}$ for all $\mathbf{p}\in \Omega_p$, $k=1,\ldots,M$. Summarizing this relation for all experimental settings, i.e., $\mathbf{S}(\mathbf{Y}(\mathbf {p}),\mathbf{p})=(\boldsymbol{S}^{(1)}(\mathbf{y}^{(1)}(\mathbf {p}),\mathbf{p}),\ldots,\boldsymbol{S}^{(M)}(\mathbf{y}^{(M)}(\mathbf {p}),\mathbf{p}))=\mathbf{0}$, we obtain $\mathrm{D}_p\mathbf{Y}$ by solving the linear system 
\begin{align*}
 \partial_y \mathbf{S}(\mathbf{Y}(\mathbf{p}),\mathbf{p}) \,\mathrm{D}_p \mathbf{Y}(\mathbf{p}) = -\partial_p \mathbf{S}(\mathbf{Y}(\mathbf{p}),\mathbf{p})
\end{align*}
The system matrix is of diagonal block form, each block $\partial_y \boldsymbol{S}^{(i)}(\mathbf{y}^{(i)}(\mathbf {p}),\mathbf{p})$ is already assembled when solving the respective fiber spinning boundary value problem. The righthand-side $\partial_p \mathbf{S}$ contains the partial derivatives with respect to $\mathbf{p}$ which we compute by means of differentiation with complex variables.

The performance of the trust region method is significantly better when both parameters $(n,\kappa)=\mathbf{p}$ are of similar magnitude. This is guaranteed by our introduction of $\kappa$ as substitute for the original $K$, $\kappa=\mathrm{ln}(K)$, cf., \eqref{eq:CarreauLikeLawK}. In the following, in the absence of more precise information, we use $K_\circ=1$~Pa as reference, cf., Remark~\ref{rem:K}, hence $10^3 < K< 10^9$ can roughly be expected. This provides a parameter domain $\Omega_p=[0,1]\times [7,20]$ for $(n,\kappa)$. And, no further scalings become necessary in the optimization scheme.

\section{Application and Results}\label{sec:results}
In the following we demonstrate our simulation-based parameter identification procedure for the elongational viscosity in high-speed spinning. As test setting we consider fiber spinning of the polymer PMMA7N in a process, where the nozzle has a length of $10^{-2}$~m and a diameter of $10^{-3}$~m. The spin-line is approximately 0.5~m long. Measurement sets of the fiber diameter are given with respect to different stencil velocities, i.e., $\{0.27, 0.53, 1.06\}$~mm/s, and take-up pressures, i.e., $\{1, 2, 3\}$~bar. For details to spinning process, measurements and data set we refer to \cite{bier2022novel}. For a complete overview of process and material parameters see Table \ref{tab:parameters}.

\subsection{Test setting: Data and solver}
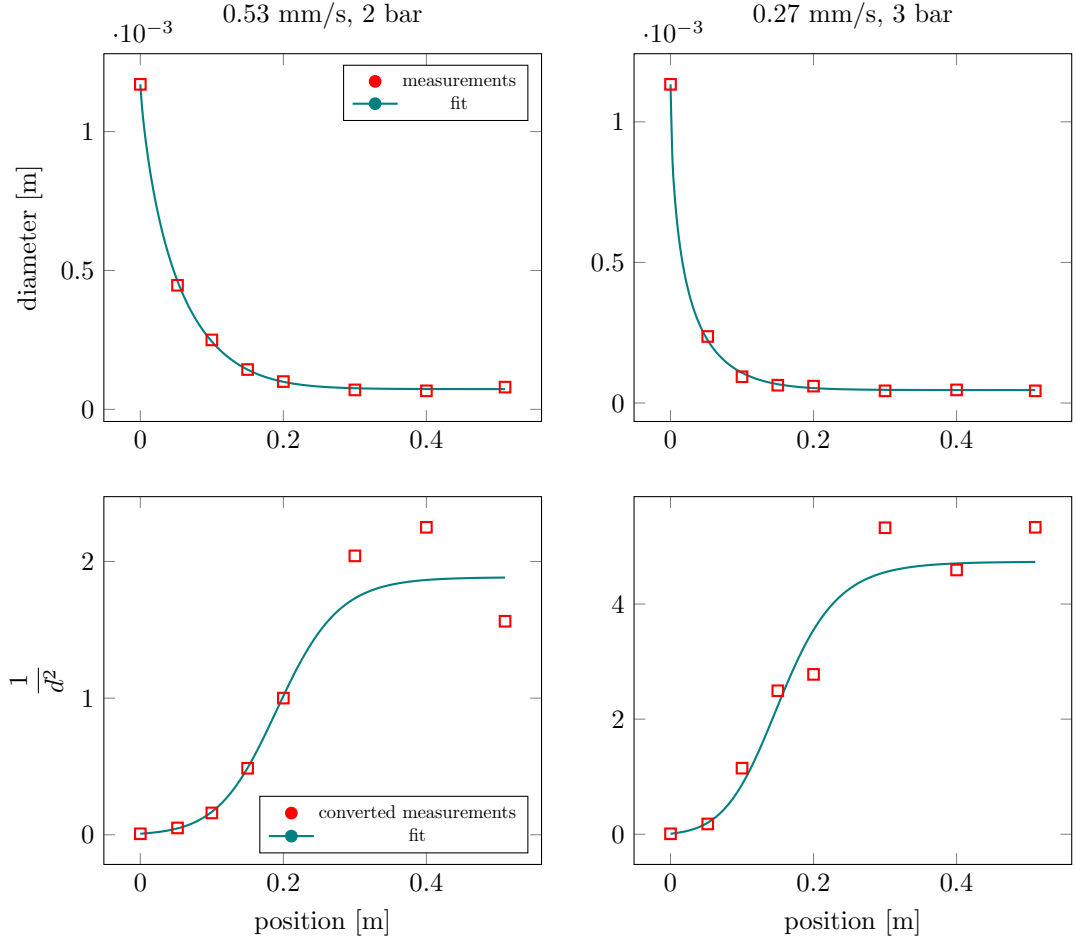
\begin{figure}[tb]
\begin{tikzpicture}
\begin{groupplot}[group style={group size=2 by 2, horizontal sep=35}, width=0.925*8cm, height=0.925*7cm, legend style={nodes={scale=0.7, transform shape}}, 
        legend image post style={mark=*},
        ]
\nextgroupplot[legend pos = north east, ylabel={diameter [m]},title={0.53 mm/s, 2 bar}]

\addplot[color=red,only marks,thick, mark=square] table[col sep=comma,header=false,x index=9,y index=4] {fiberDiameter.csv};
\addlegendentry{measurements}

\addplot[color=teal,smooth,thick] table[col sep=comma,header=false,x index=9,y index=4] {fittedDiameterProfile.csv};
\addlegendentry{fit}

\nextgroupplot[legend pos = south east,  ,title={0.27 mm/s,  3 bar}]

\addplot[color=red,only marks,thick, mark=square] table[col sep=comma,header=false,x index=9,y index=8] {fiberDiameter.csv};

\addplot[color=teal,smooth,thick] table[col sep=comma,header=false,x index=9,y index=8] {fittedDiameterprofile.csv};


\nextgroupplot[legend pos = south east,  xlabel={position [m]}, ylabel={$\dfrac{1}{d^2}$}]

\addplot[color=red,only marks, mark=square, thick] table[col sep=comma,header=false,x index=9,y index=4] {scaledConvertedMeasurements.csv};
\addlegendentry{converted measurements}

\addplot[color=teal,smooth,thick] table[col sep=comma,header=false,x index=9,y index=4] {fittedConvertedDiameterProfile.csv};
\addlegendentry{fit}


\nextgroupplot[ legend pos = north east, 
        xlabel={position [m]} ]

\addplot[color=red,only marks,thick, mark=square] table[col sep=comma,header=false,x index=9,y index=8] {scaledConvertedMeasurements.csv};

\addplot[color=teal,smooth,thick] table[col sep=comma,header=false,x index=9,y index=8] {fittedConvertedDiameterProfile.csv};        

\end{groupplot}

\end{tikzpicture}
\caption{Data of measured fiber diameters (top) and converted (velocity-like) counterparts (bottom) with associated fitted profiles along the spin-line for two different experimental settings (left and right)}
\label{fig:FittedMeasurements}
\end{figure}
The underlying data on the measured diameters refer to nine experiments with different stencil velocities and take-up pressures. The experiments cover ranges of high and low strain rates. We observe that the measurements are more susceptible to disturbances at lower strain rates and near the take-up point due to the used aspirator. As example, Fig.~\ref{fig:FittedMeasurements} shows the measured diameters along the spin-line for two spun fibers. The disturbances become more evident in the converted (dimensionless) counterparts $d^{-2}$ that can be interpreted as velocity. Our proposed smoothing strategy \eqref{eq:ansatzFunc} provides a continuous diameter profile that captures the behavior near the nozzle very well and suitably balances the noisy data near the take-up point.

Concerning our model-based simulation framework, we refer to Appendix \ref{sec:ModelClosure} for the specification of the used constitutive laws, closure relations and airflow model. The chosen referential values for the non-dimensionalization are summarized in Table \ref{tab:referenceValues}. Whereas, for the boundary conditions, the inlet information $u_{in}$, $T_{in}$ is well known from the experimental settings, we adapt the outlet information $d_{out}$ from the noisy measurements caused by the aspirator. In particular, we set $d_{out}=d^{(k)}_{\mathrm{fit}}(L)$. We note that the overall simulated fiber behavior turns out to be robust to small perturbations of the boundary conditions for $d_{out}$, although the identified values of the material parameters may slightly change.

We realize our parameter identification in MATLAB. We use the core of the MATLAB routine \texttt{lsqnonlin} which implements a subspace trust region method based on \cite{coleman1996interior} and use the predefined termination criteria. The occurring boundary value problems are solved by means of a collocation-continuation scheme, where the MATLAB rountine \texttt{bvp4c} is used for the collocation (Lobatto IIIa), (cf.,  Appendix \ref{sec: NumericalBVPs}). The evaluation of the cost function $J$ \eqref{eq:opt:disc} is performed on regular grids $\Delta_o^{(k)}$ with $n_o^{(k)}=20$.
All simulations are performed on a machine with Intel(R) Core(TM) i7-10510U CPU @ 1.80GHz and 32 GB Ram using Matlab Version R2021a.

\subsection{Identified material parameters and fiber behavior}
For the given data set we find the material parameters $\mathbf{p}_\mathrm{opt}=(n,\kappa)_\mathrm{opt}=(0.8,11.71)$ as optimal solution of \eqref{eq:opt:disc}, implying ${K_\mathrm{opt}=1.22\cdot 10^{5}}$~Pa.
The associated diameter profiles and velocity-like counterparts are visualized in 
Fig.~\ref{fig:OptimaPlots} and Fig.~\ref{fig:OptimaPlots_New} for several experimental set-ups. 
The parameters $\mathbf{p}_\mathrm{opt}$ provide a Carreau-like model with pronounced non-Newtonian behavior. 
The differences to the Newtonian behavior ($n=1$) are marginal for small strain rates, but become clearly visible in the experimental set-ups with higher take-up pressures which cause significant strain rates.
The identified material law captures the smoothed measurements in the settings with higher stencil velocities very well. For lower stencil velocities we observe an overestimate in the velocity-like profiles. In comparison, the Newtonian simulations overestimate the velocity-like profiles in all settings; this is accompanied by a general underestimation of the diameter profiles. Summing up, the overall approximation quality of the simulation results obtained with the optimized Carreau-like material law is obviously better than the one with the Newtonian material law ($n=1$). Our results suggest a non-Newtonian behavior of the considered polymer PMMA7N which confirms the results of the more heuristic approach proposed in \cite{bier2022novel}.

\begin{figure}[tb]
\begin{center}
\begin{tikzpicture}
\begin{groupplot}[group style={group size=2 by 3, horizontal sep=35, vertical sep=40}, width=7.5cm, height=5.5cm, legend style={nodes={scale=0.6, transform shape}}, 
legend image post style={mark=*},
]

\nextgroupplot[title= {1.06 mm/s, 2 bar},
               ylabel={diameter [m]},] 
\addplot[color=red,only marks,thick, mark=square] table[col sep=comma,header=false,x index=9,y index=1] {fiberDiameter.csv};
\addlegendentry{measured}

\addplot[color=red,smooth, dashed] table[col sep=comma,header=false,x index=9,y index=1] {fittedDiameterProfile.csv};
\addlegendentry{fit}

\addplot[color=black,smooth] table[col sep=comma,header=false,x index=9,y index=1] {newtonDiams.csv};
\addlegendentry{simulated, $n=1$}

\addplot[color=green,smooth,thick] table[col sep=comma,header=false,x index=9,y index=1] {optDiams.csv};
\addlegendentry{simulated, $\mathbf{p}_\mathrm{opt}$}

\nextgroupplot[ title= {1.06 mm/s, 3 bar}
               ] 
\addplot[color=red,only marks,thick, mark=square] table[col sep=comma,header=false,x index=9,y index=2] {fiberDiameter.csv};

\addplot[color=red,smooth, dashed] table[col sep=comma,header=false,x index=9,y index=2] {fittedDiameterProfile.csv};

\addplot[color=black,smooth] table[col sep=comma,header=false,x index=9,y index=2] {newtonDiams.csv};

\addplot[color=green,smooth,thick] table[col sep=comma,header=false,x index=9,y index=2] {optDiams.csv};

\nextgroupplot[title= {0.53 mm/s, 2 bar},
                legend pos = north east,  ylabel={diameter [m]},]
\addplot[color=red,only marks,thick, mark=square] table[col sep=comma,header=false,x index=9,y index=4] {fiberDiameter.csv};

\addplot[color=red,smooth, dashed] table[col sep=comma,header=false,x index=9,y index=4] {fittedDiameterProfile.csv};

\addplot[color=black,smooth] table[col sep=comma,header=false,x index=9,y index=4] {newtonDiams.csv};

\addplot[color=green,smooth,thick] table[col sep=comma,header=false,x index=9,y index=4] {optDiams.csv};

\nextgroupplot[title= {0.53 mm/s, 3 bar},
                legend pos = north east, 
       ]
\addplot[color=red,only marks,thick, mark=square] table[col sep=comma,header=false,x index=9,y index=5] {fiberDiameter.csv};

\addplot[color=red,smooth, dashed] table[col sep=comma,header=false,x index=9,y index=5] {fittedDiameterProfile.csv};

\addplot[color=black,smooth] table[col sep=comma,header=false,x index=9,y index=5] {newtonDiams.csv};

\addplot[color=green,smooth,thick] table[col sep=comma,header=false,x index=9,y index=5] {optDiams.csv};

\nextgroupplot[title= {0.27 mm/s, 2 bar},
                legend pos = north east,  xlabel={position [m]}, ylabel={diameter [m]},]
\addplot[color=red,only marks,thick, mark=square] table[col sep=comma,header=false,x index=9,y index=7] {fiberDiameter.csv};

\addplot[color=red,smooth, dashed] table[col sep=comma,header=false,x index=9,y index=7] {fittedDiameterProfile.csv};

\addplot[color=black,smooth] table[col sep=comma,header=false,x index=9,y index=7] {newtonDiams.csv};

\addplot[color=green,smooth,thick] table[col sep=comma,header=false,x index=9,y index=7] {optDiams.csv};

\nextgroupplot[title= {0.27 mm/s, 3 bar}, 
                legend pos = north east,  xlabel={position [m]}, ]
\addplot[color=red,only marks,thick, mark=square] table[col sep=comma,header=false,x index=9,y index=8] {fiberDiameter.csv};

\addplot[color=red,smooth, dashed] table[col sep=comma,header=false,x index=9,y index=8] {fittedDiameterProfile.csv};

\addplot[color=black,smooth] table[col sep=comma,header=false,x index=9,y index=8] {newtonDiams.csv};

\addplot[color=green,smooth,thick] table[col sep=comma,header=false,x index=9,y index=8] {optDiams.csv};

\end{groupplot}
\end{tikzpicture}
\end{center}
\caption{Diameter profiles over spin-line in different experimental set-ups for PMMA7N: measured data and associated fit vs. simulated results with $\mathbf{p}_\mathrm{opt}$ as well as $n=1$ (Newtonian)}
\label{fig:OptimaPlots}
\end{figure}
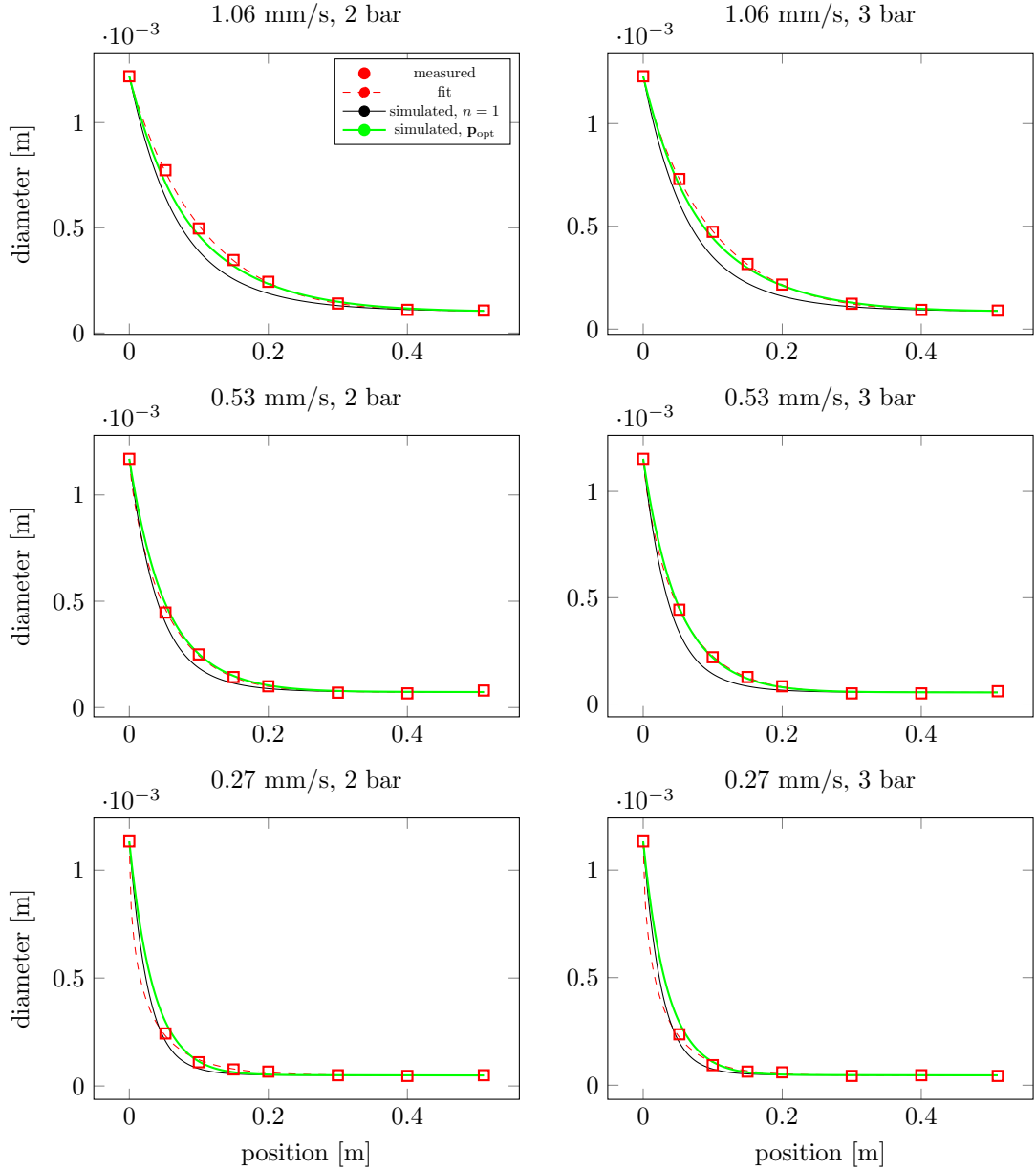
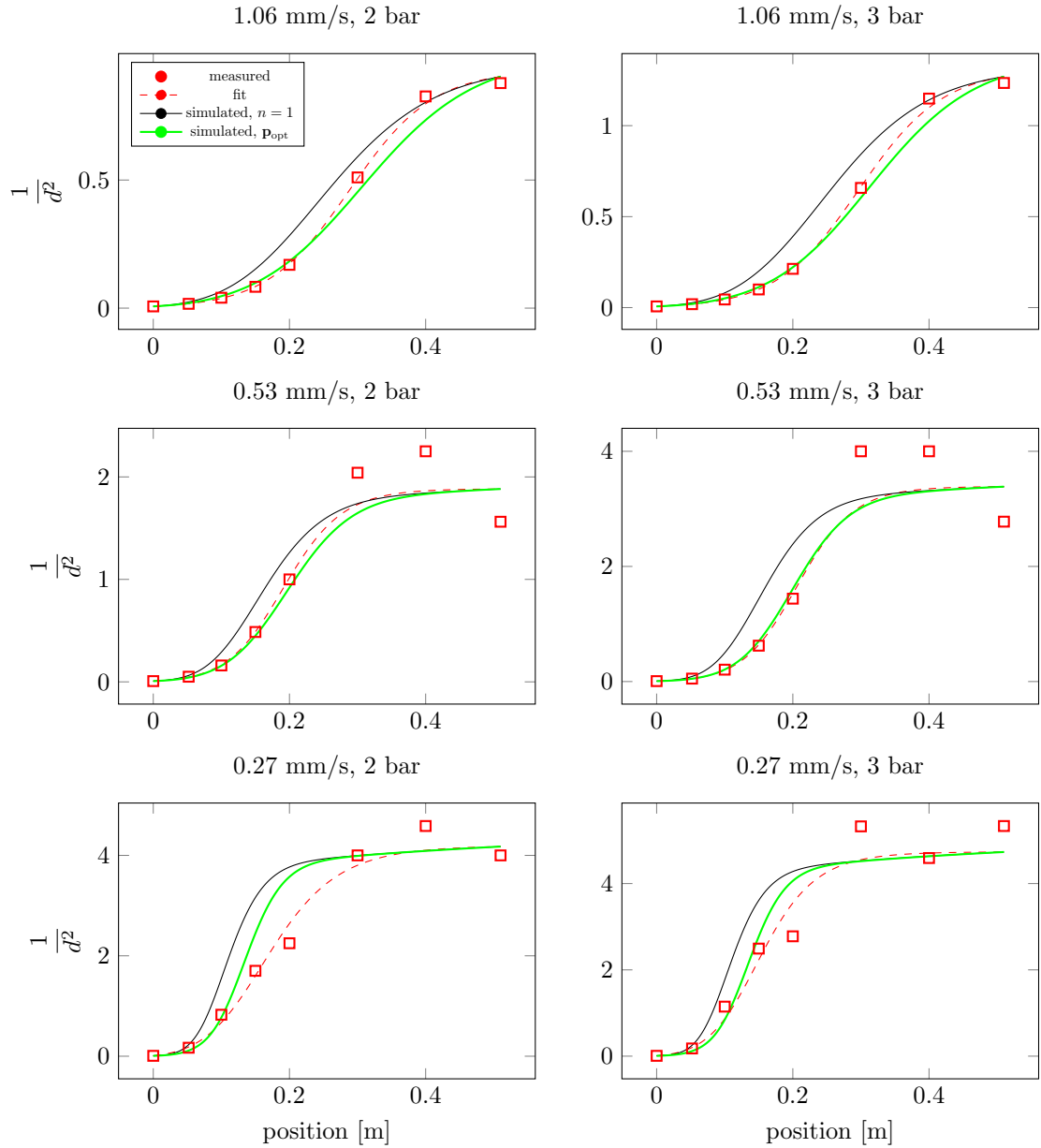
\begin{figure}[tb]
\begin{center}
\begin{tikzpicture}
\begin{groupplot}[group style={group size=2 by 3, horizontal sep=35, vertical sep=40}, width=7.5cm, height=5.5cm, legend style={nodes={scale=0.6, transform shape}}, 
        legend image post style={mark=*},
        ]

\nextgroupplot[title= {1.06 mm/s, 2 bar},
               ylabel={$\dfrac{1}{d^2}$},
               legend pos = {north west}] 
\addplot[color=red,only marks,thick, mark=square] table[col sep=comma,header=false,x index=9,y index=1] {scaledConvertedMeasurements.csv};
\addlegendentry{measured}

\addplot[color=red,smooth, dashed] table[col sep=comma,header=false,x index=9,y index=1] {fittedConvertedDiameterProfile.csv};
\addlegendentry{fit}

\addplot[color=black,smooth] table[col sep=comma,header=false,x index=9,y index=1] {convertedNewtonDiams.csv};
\addlegendentry{simulated, $n=1$}

\addplot[color=green,smooth,thick] table[col sep=comma,header=false,x index=9,y index=1] {convertedOptDiams.csv};
\addlegendentry{simulated, $\mathbf{p}_\mathrm{opt}$}

\nextgroupplot[ title= {1.06 mm/s, 3 bar}
               ] 
\addplot[color=red,only marks,thick, mark=square] table[col sep=comma,header=false,x index=9,y index=2] {scaledConvertedMeasurements.csv};

\addplot[color=red,smooth, dashed] table[col sep=comma,header=false,x index=9,y index=2] {fittedConvertedDiameterProfile.csv};

\addplot[color=black,smooth] table[col sep=comma,header=false,x index=9,y index=2] {convertedNewtonDiams.csv};

\addplot[color=green,smooth,thick] table[col sep=comma,header=false,x index=9,y index=2] {convertedOptDiams.csv};

\nextgroupplot[title= {0.53 mm/s, 2 bar},
                legend pos = north east,  ylabel={$\dfrac{1}{d^2}$ },]
\addplot[color=red,only marks,thick, mark=square] table[col sep=comma,header=false,x index=9,y index=4] {scaledConvertedMeasurements.csv};

\addplot[color=red,smooth, dashed] table[col sep=comma,header=false,x index=9,y index=4] {fittedConvertedDiameterProfile.csv};

\addplot[color=black,smooth] table[col sep=comma,header=false,x index=9,y index=4] {convertedNewtonDiams.csv};

\addplot[color=green,smooth,thick] table[col sep=comma,header=false,x index=9,y index=4] {convertedOptDiams.csv};

\nextgroupplot[title= {0.53 mm/s, 3 bar},
                legend pos = north east, 
       ]
\addplot[color=red,only marks,thick, mark=square] table[col sep=comma,header=false,x index=9,y index=5] {scaledConvertedMeasurements.csv};

\addplot[color=red,smooth, dashed] table[col sep=comma,header=false,x index=9,y index=5] {fittedConvertedDiameterProfile.csv};

\addplot[color=black,smooth] table[col sep=comma,header=false,x index=9,y index=5] {convertedNewtonDiams.csv};

\addplot[color=green,smooth,thick] table[col sep=comma,header=false,x index=9,y index=5] {convertedOptDiams.csv};

\nextgroupplot[title= {0.27 mm/s, 2 bar},
                legend pos = north east,  xlabel={position [m]}, ylabel={$\dfrac{1}{d^2}$},]
\addplot[color=red,only marks,thick, mark=square] table[col sep=comma,header=false,x index=9,y index=7] {scaledConvertedMeasurements.csv};

\addplot[color=red,smooth, dashed] table[col sep=comma,header=false,x index=9,y index=7] {fittedconvertedDiameterProfile.csv};

\addplot[color=black,smooth] table[col sep=comma,header=false,x index=9,y index=7] {convertedNewtonDiams.csv};

\addplot[color=green,smooth,thick] table[col sep=comma,header=false,x index=9,y index=7] {convertedOptDiams.csv};

\nextgroupplot[title= {0.27 mm/s, 3 bar}, 
                legend pos = north east,  xlabel={position [m]}, ]
\addplot[color=red,only marks,thick, mark=square] table[col sep=comma,header=false,x index=9,y index=8] {scaledConvertedMeasurements.csv};

\addplot[color=red,smooth, dashed] table[col sep=comma,header=false,x index=9,y index=8] {fittedConvertedDiameterProfile.csv};

\addplot[color=black,smooth] table[col sep=comma,header=false,x index=9,y index=8] {convertedNewtonDiams.csv};

\addplot[color=green,smooth,thick] table[col sep=comma,header=false,x index=9,y index=8] {convertedOptDiams.csv};

\end{groupplot}
\end{tikzpicture}
\end{center}
\caption{Velocity-like profiles over spin-line in different experimental set-ups for PMMA7N: converted measured data and associated fit vs. simulated results with $\mathbf{p}_\mathrm{opt}$ as well as $n=1$ (Newtonian)}
\label{fig:OptimaPlots_New}
\end{figure}

\begin{figure}[tb]
\includegraphics[width=\textwidth]{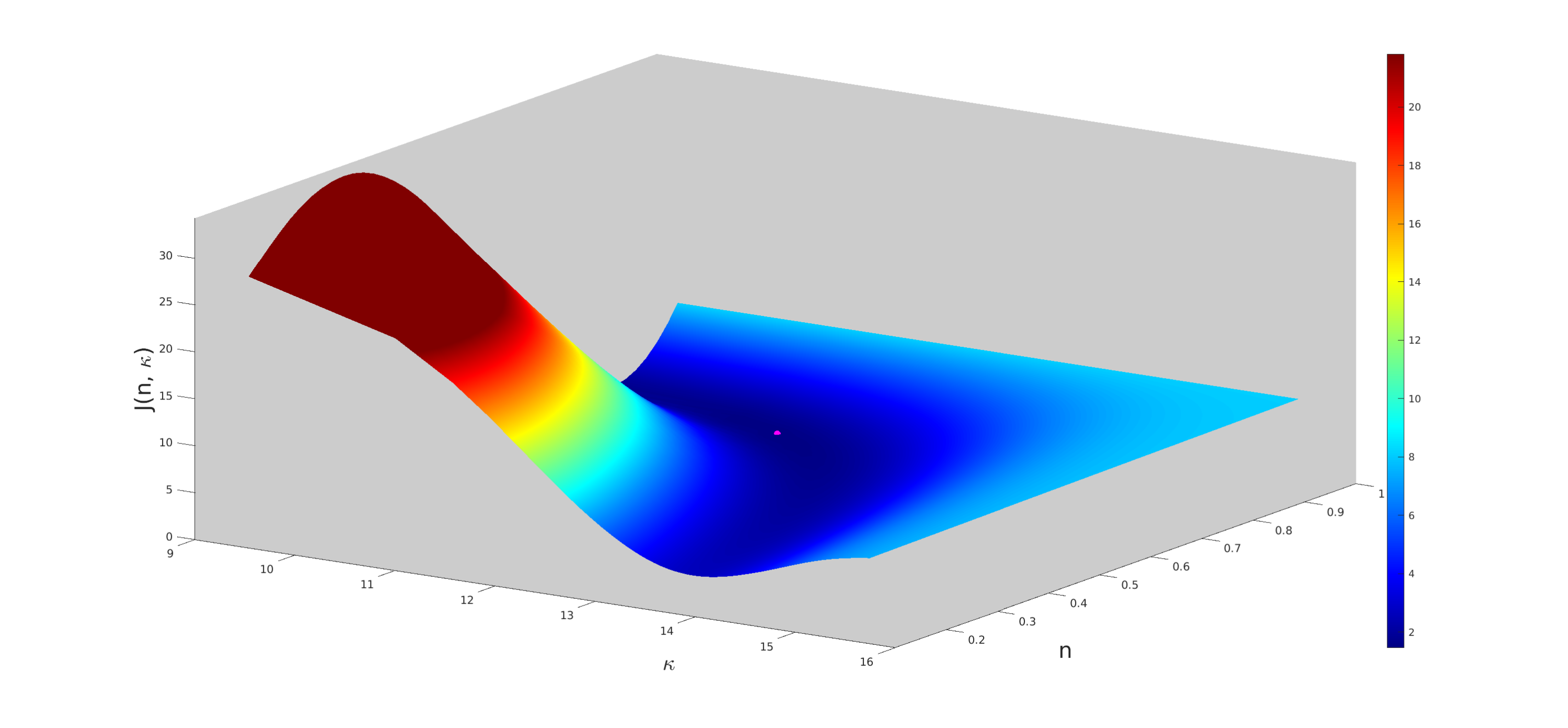}
    \caption{Cost function $J$ for $n\in [0,1]$ and $\kappa \in [9,16]$}
    \label{fig:costFunction}
\end{figure}

From a mathematical point of view we note that $\mathbf{p}_\mathrm{opt}$ is a global minimizer of the cost function $J$. In the considered application, the parameter identification is a two-parametric optimization problem. The underlying low-dimensional parameter space allows the resolution of $J$, for a visualization of the cost function see Fig.~\ref{fig:costFunction}. For $n\rightarrow 1$ or large $\kappa$, implying $K\rightarrow \infty$, the cost function $J$ reaches a plateau. This comes from the fact that the Carreau-like material law reduces to a strain rate-independent Vogel-Fulcher-Tammann law for a Newtonian fluid, i.e., $\mu_e(T,\dot\epsilon;(n,\kappa)) \xrightarrow{{n\rightarrow 1}} \mu_{e,0}(T)$ as well as $\mu_e(T,\dot\epsilon;(n,\kappa)) \xrightarrow{{\kappa\rightarrow\infty}} \mu_{e,0}(T)$. For $n<1$ and $\kappa$ small enough, the effect of the strain rate-dependent term in $\mu_{e}$ becomes visible, yielding smaller values of $J$. However, if $n$ and $\kappa$ are chosen too small, the non-Newtonian effect is too strong and the associated simulated fiber diameters strongly overestimate the measurements. Thus, the values of the cost function rise again and become even higher than in the Newtonian plateau. Furthermore, it is possible that the simulations fail: the boundary value problems are not solvable any more as a consequence of the exploding strain rates at the fiber end due to the unsuitable parameterization. Figure~\ref{fig:costFunction} clearly shows that there is a parameter range, for which the associated Carreau-like viscosity law outperforms the Vogel-Fulcher-Tammann one of a Newtonian fluid. 

With respect to the approach in [1], our results are qualitatively comparable, but of course not quantitatively. This is due to the fact that in [1] no numerical simulations of a non-Newtonian material were performed. Moreover, for density and specific heat capacity, information from a constant model was considered. The determined parameters $n=0.22$ and $K=9.4\cdot 10^{5}$ Pa are significantly lower in $n$ and higher in $K$ than ours, implying $\kappa=13.75$. Be that as it may, this parameter combination does not provide a bad approximation, as can be seen from Fig.~\ref{fig:costFunction}.

\subsection{Solver performance}
Any evaluation of the cost function requires the solution of multiple boundary value problems. Therefore, a detailed resolution of the cost function (as in Fig.~\ref{fig:costFunction}) is computationally expensive -- and impractical in the general case of underlying higher-dimensional parameter spaces. The proposed gradient-based optimization procedure requires only the computation of a few function values to identify the material parameters. It is therefore also applicable to other higher-dimensional parameter identification problems. 

\begin{table}[t]
\caption{Performance of optimization scheme for different initial guess $\mathbf{p}_\mathrm{init}\in \Omega_p$}
\begin{tabular}{|r r||r r |r|r|}
\hline
\multicolumn{2}{|c||}{Initial $\mathbf{p}_{\mathrm{init}}=(n,\kappa)_\mathrm{init}$} & \multicolumn{2}{|c|}{Identified $\mathbf{p}_{\mathrm{opt}}=(n,\kappa)_\mathrm{opt}$} & Iteration count & CPU [s] \\ \hline
0.5 & 9.21 &0.7996&11.7105&  8         & 229                 \\
0.7 & 9.21 &0.7996&11.7108&   4          & 64                  \\
0.9 & 9.21 &0.7994&11.7142&   12          & 86                  \\
0.3 & 11.51 &0.7997&11.7092& 9         &  178            \\
0.5 & 11.51 &0.7996&11.7106&  10          & 120                  \\
0.7 & 11.51 &0.7995&11.7120&  10          & 79                  \\
0.9 & 11.51 &0.7994&11.7141&  9          & 75                  \\
0.3 & 13.82 &0.7993&11.7150& 12          & 100                  \\
0.5 & 13.82 &0.7997&11.7116& 12          & 94                  \\
0.7 & 13.82 &0.7994& 11.7140 &9          & 82                 \\
0.9 & 13.82 &0.7995&11.7121& 10          & 87                  \\ 
\hline
\end{tabular}
\label{tab:perfResults}
\end{table}

\IncMargin{1em}
\begin{algorithm}[b]
\SetKwData{Left}{left}\SetKwData{This}{this}\SetKwData{Up}{up}
\SetKwFunction{Union}{Union}\SetKwFunction{FindCompress}{FindCompress}
\SetKw{Continue}{continue}
\SetKwInOut{Input}{input}\SetKwInOut{Output}{output}
\Input{$J:$ objective function \\ 
$n_\mathrm{init}$: initialization of $n$, $n_\mathrm{init}<1$ \\
$\kappa_l$, $\kappa_u$: lower and upper bound for $\kappa$}
\Output{starting point $\mathbf{p}_\mathrm{init}=(n,\kappa)$}
\BlankLine

$\kappa \leftarrow \kappa_u$ \;
$J_1 \leftarrow J(1,\kappa)$ \;
$J_2 \leftarrow J(n_\mathrm{init},\kappa)$ \;
$relDif \leftarrow \left | \dfrac{J_1 - J_2}{J_1} \right| $ \;
\While{($relDif < \frac{1}{10}$) and ($\kappa> \kappa_l$)}{
 $\kappa \leftarrow \kappa-1 $\;
$J_2 \leftarrow$ $J(n_\mathrm{init},\kappa)$ \;
 $relDif \leftarrow \left| \dfrac{J_1 - J_2}{J_1} \right |$ \;
}
$\mathbf{p}_\mathrm{init} \leftarrow (n_\mathrm{init},\kappa)$ \;
\caption{Starting point heuristic}\label{algo:startHeuristic}
\end{algorithm}\DecMargin{1em}

The performance of the optimization scheme strongly depends on balanced parameters of similar magnitude and on a suitable initial guess / starting point $\mathbf{p}_\mathrm{init}$. With respect to the Carreau-like model, the balance of the material parameters is achieved here by considering $\kappa=\mathrm{ln}(K)$ instead of $K$. Starting in the neighborhood of the minimizer, the trust region method quickly converges to it within a few iterations, see Table~\ref{tab:perfResults}. Each iteration requires the evaluation of cost function as well as first and second derivatives, cf., Sec.~\ref{subsec:opt-alg}, and takes a CPU time of about 9 seconds. Note that a speed up could be achieved by parallelization, as the underlying $M$ boundary value problems are solved sequentially here. However, in case of an unsuitable initial guess, the optimization method might fail. We avoid trial-and-error-strategies and instead propose a heuristic to determine an appropriate $\mathbf{p}_\mathrm{init}$. Without knowledge about the exact profile of the cost function, the consideration of some limits allows for an efficient heuristic procedure. The Newtonian limit ($n=1$ or $K=\exp(\kappa)\rightarrow \infty$) can be excluded as starting point as there are zero directional derivatives with respect to $\kappa$ in case of $n=1$ as well as almost zero directional derivatives with respect to $\kappa$ and $n$ in case of large $\kappa$, causing a stagnation of the optimization scheme. However, as the Carreau-like viscosity model simplifies to the Vogel-Fulcher-Tammann one for a Newtonian fluid as $n\rightarrow 1$ and $K\rightarrow \infty$, it is promising to reverse the transition and look for a starting point by slowly adding the non-Newtonian effects. Hence, we fix $n$ to a moderate value,  $n=n_\mathrm{init}<1$, and use the Newtonian results for a large $\kappa$ as reference. Then, we successively decrease $\kappa$ by a factor 1 and evaluate the cost function, until we observe significant changes in $J$. The resulting parameter tuple is then used a initial guess for the trust region method. The additional computational effort is rather low. The boundary values problems occurring in the cost function evaluations are almost identical and thus efficiently solvable, if the Newtonian results are used for initialization. The start point heuristic is depicted in Algorithm~\ref{algo:startHeuristic}. We use here $n_\mathrm{init}=0.5$ and $\kappa\in [\kappa_l,\kappa_b]=[7,20]$.
Note, while the optimization scheme is directly applicable to other, also higher-dimensional, parameter identification problems, the proposed start point heuristic is problem-tailored and must be adapted if the material law or its parameterization is changed.

\section{Conclusion}
This work presents a simulation-based procedure for the identification of material parameters in fiber spinning processes. The study acts as a proof of concept. Presupposing a Carreau-type material law,
the elongational viscosity of PMMA is determined in a high-speed setup on top of fiber diameter measurements. The gradient-based optimization shows to be robust and efficient. The presented framework is straightforward applicable to generalized Newtonian fluids with parametrically described material law.

\appendix
\section{Model Closure}\label{sec:ModelClosure}
The appendix provides the used models for the aerodynamic forces $f_{air}$ and the heat transfer coefficient $\alpha$, as well as the models for the density $\rho$ and the specific heat capacity $c_p$ of the considered material. Moreover, Table~\ref{tab:parameters} and Table~\ref{tab:referenceValues} give an overview of the process and physical parameters as well as the reference values used in the non-dimensionalization.

To distinguish between the dimensional and non-dimensional quantities, we proceed like in the beginning of Sec.~\ref{sec:FiberModel}, where the dimensionless counterpart to the dimensional quantity $x$ is denoted by   $\tilde x$.
The aerodynamic forces $f_{air}$ are described by the air drag model $\mathbf{F}$ from \cite{marheineke2011modeling}.
The dimensional formulation for a general flow situation reads
$$ \mathbf{f}_{air}(\mathbf{\boldsymbol{\tau}}, \mathbf{v}_{rel}, \nu_*, \rho_*, d) = \frac{ \rho_\star \nu^2_*}{d}\, \mathbf{F}(\boldsymbol{\tau}, \frac{d}{\nu_*}\,\mathbf{v}_{rel})$$
with (normalized) fiber tangent $\boldsymbol{\tau}$, relative velocity between air and fiber $\mathbf{v}_{rel}= \mathbf v_* - \mathbf v$ and fiber diameter $d$. The aerodynamic quantities are marked with $_\star$, i.e., air density $\rho_*$ and kinematic viscosity $\nu_*$. 
The dimensionless air drag function $\mathbf{F}$ is modeled in terms of the tangential and normal components of the relative velocity, i.e.,
\begin{align*}
    \mathbf{F}(\boldsymbol{\tau}, \mathbf{v}) = v_n r_n^\delta(v_n) \mathbf{n} + v_{\tau} r_\tau^\delta(v_n)\boldsymbol{\tau},
\end{align*}
with $v_\tau = \mathbf{v} \cdot \boldsymbol{\tau}$,  $v_n = \mathbf v \cdot \mathbf n$ and $\mathbf{n}=\mathbf{v}-v_\tau \boldsymbol{\tau}/\|\mathbf{v}-v_\tau\boldsymbol{\tau}\|$.
The resistance coefficients $r_n^\delta$ and $r_\tau^\delta$ are regularized with respect to the slenderness ratio $\delta$.
The non-dimensional formulation 
\begin{align*}\tilde{\mathbf{f}}_{air}(\boldsymbol{\tau}, \mathbf{\tilde v}_{rel}, \tilde \nu_*, \tilde \rho_*, \tilde d) = \frac{\mathrm{A}_*}{\mathrm{Re}_*^2}\frac{\tilde{\rho}_\star\tilde{\nu}_*^2}{\tilde{d}}\,\mathbf{F}(\boldsymbol{\tau}, \mathrm{Re}_* \frac{\tilde{d}}{\tilde{\nu}_*}\,\tilde{\mathbf{v}}_{rel}),\qquad \mathrm{A}_*= \frac{\rho_{*,\circ}  d_\circ u_\circ^2}{f_\circ}, \quad \mathrm{Re_*}= \frac{ d_\circ u_\circ}{\nu_{*,\circ}}
\end{align*}
relies on the characteristic air-associated Reynolds number $\mathrm{Re}_\star $ and the air drag number $\mathrm{A}_\star$ with $f_\circ=Q_\circ/(L_\circ u_\circ)$.

The heat transfer coefficient $\alpha$ depends additionally on the specific heat capacity $c_{p,\star}$ and the heat conductivity $\lambda_\star$ of the air. It is modeled in terms of the  Nusselt function $\mathcal{N}$,  cf., \cite{wieland2019,sucker1976steady}. The non-dimensional variant relies  on the characteristic air-associated Reynolds number $\mathrm{Re}_\star $, the characteristic Nusselt number $\mathrm{Nu}_\star$ and the Prandtl number $\mathrm{Pr}_\star$,
\begin{align*}
\alpha(\boldsymbol{\tau}, \mathbf{v}_{rel},\nu_*, \rho_*,  c_{p,*}, \lambda_*,  d) &= \frac{\lambda_*}{d} \mathcal{N}(\dfrac{d}{\nu_*} {\mathbf{v}}_{rel} \cdot \boldsymbol{\tau},\dfrac{d}{\nu_*}\Vert {\mathbf{v}}_{rel}\Vert_2,  \frac{ c_{p,*} \rho_\star \nu_*}{\lambda_*} ),
\\
\tilde{\alpha}(\boldsymbol{\tau}, \mathbf{\tilde{v}}_{rel},  \tilde \nu_*,  \tilde \rho_*, \tilde c_{p,*}, \tilde \lambda_*, \tilde d) &= \frac{1}{\mathrm{Nu}_*}\frac{\tilde \lambda_*}{\tilde d} \mathcal{N}(\mathrm{Re}_*\dfrac{\tilde d} {\tilde \nu_*} \tilde{\mathbf{v}}_{rel} \cdot \boldsymbol{\tau},\mathrm{Re}_*\dfrac{\tilde d} {\tilde \nu_*}\Vert \tilde{\mathbf{v}}_{rel}\Vert_2, \mathrm{Pr}_* \frac{\tilde c_{p,*} \tilde{\rho}_\star \tilde \nu_*}{\tilde \lambda_*} ),\\
&\qquad \mathrm{Nu}_* = \frac{\alpha_\circ d_\circ}{\lambda_{*,\circ}}, \quad \mathrm{Pr}_* = \frac{c_{p,\star,\circ} \rho_{\star,\circ} \nu_{\star.\circ}}{\lambda_{*,\circ}}
\end{align*}

In the spinning process under consideration, the air is at rest ($\mathbf{v}_\star =\mathbf{0}$) and the fiber is straight, so $\mathbf{v}_{rel}=-u\boldsymbol{\tau}$. Hence, the air drag model simplifies to a scalar-valued function with Stokes resistance coefficient $r^\delta_{\tau,S}$, i.e.,
$$f_{air}(u,d,s)=
 -u\mu_*(s)\,r_{\tau,S}^\delta, \qquad r_{\tau,S}^\delta = \frac{2\pi}{\ln(\frac{4}{\delta})} + \frac{{\pi}/{2}}{\ln^2(\frac{4}{\delta})}
 $$
 We use here $\delta = 10^{-3}$. The heat transfer coefficient becomes
 $$\alpha(u,d,s)=
  \frac{\lambda_*(s)}{d} \mathcal{N}(-\dfrac{du}{\nu_*(s)},\dfrac{du}{\nu_*(s)},  \frac{ c_{p,*} \rho_\star \nu_*}{\lambda_*}(s) )
 $$
 Moreover, mass density, specific heat capacity, kinematic viscosity and heat conductivity of the air are taken here as constant. The exact values can be found in Table \ref{tab:parameters}.

\begin{table}[t] 
\caption{Process and physical parameters}
\setlength{\tabcolsep}{1em}
\begin{tabular}{|llll|l}
\cline{1-4}
Description                        & Symbol     & Value     & Unit                        &  \\ \cline{1-4}
Fiber length                       & $L$          & 0.51       & m                           &  \\
Temperature at nozzle              & $T_{in}$  & 513.15    & K                           &  \\
Density coefficient                & $a_\varrho$          & -0.964        & kg/($\mathrm{m}^3 \mathrm{K}$) &  \\
Density coefficient                & $b_\varrho$          & 1572.33      & kg/($\mathrm{m}^3$)   &  \\
Specific heat capacity coefficient & $a_{c_p}$          & 3.2         & J/(kg$\mathrm{K}^2$)  &  \\
Specific heat capacity coefficient & $b_{c_p}$          & 648.22       & J/(kgK)                     &  \\ 
VFT parameter & $\mu_c$ & $3.7074 \cdot 10^{-4}$ &Pa$\cdot$s\\
VFT parameter & $B$ & 3649 &K\\
VFT parameter & $T_{VFT}$ & 273.15 & K\\
Air mass density & $\rho_{*}$ & 1 & kg/$\mathrm{m}^3$\\
Air specific heat capacity coefficient & $c_{p,*}$ & 1000 &J/(kgK) \\
Air kinematic viscosity & $\nu_{*}$ & $2 \cdot 10^{-5}$ & $\mathrm{m}^2$/s\\
Air heat conductivity & $\lambda_{*}$ & 0.031 & W/(mK)\\

\cline{1-4}
\end{tabular}
\label{tab:parameters}
\end{table}%

\begin{table}[t] 
\caption{Reference values for non-dimensionalization}

\setlength{\tabcolsep}{1em}
\begin{tabular}{|llll|l}
\cline{1-4}
Description                        & Symbol     & Value     & Unit                        &  \\ \cline{1-4}
Mass flow                     & $Q_\circ$          & $3.08 \cdot 10^{-5}$       & kg/s                           &  \\
Length scale              & $L_\circ$  & $0.51$    & m                           &  \\
Velocity                & $u_\circ$          & $0.0283$        & m/s &  \\
Temperature                & $T_\circ$          & $513.15$      & K   &  \\
Mass density & $\rho_\circ$          & $1.077 \cdot 10^3$         & kg/$\mathrm{m}^3$  &  \\
Specific heat capacity coefficient & $c_{p, \circ}$          & $2.2903 \cdot 10^3$       & J/(kgK)                     &  \\ 
Dynamic viscosity & $\mu_\circ$ & $1.4865 \cdot 10^{3}$ &Pa$\cdot$s\\
Laminar heat transfer coefficient & $\alpha_\circ$ & 12.762 &W/($\mathrm{m}^2$K)\\
Air mass density & $\rho_{*, \circ}$ & 1 & kg/$\mathrm{m}^3$\\
Air specific heat capacity coefficient & $c_{p,*, \circ}$ & 1000 & J/(kgK)  \\
Air kinematic viscosity & $\nu_{*, \circ}$ & $2 \cdot 10^{-5}$ & $\mathrm{m}^2/\mathrm{s}$\\
Air heat conductivity & $\lambda_{*, \circ}$ & 0.031 & W/(mK)\\
\cline{1-4}
\end{tabular}
\label{tab:referenceValues}
\end{table}%

The material laws for the mass density $\rho$ and the 
specific heat capacity $c_p$ of the polymer are assumed to be linearly dependent on the fiber temperature $T[\mathrm{K}]$, i.e., 
\begin{align*}
g(T) = a_gT + b_g, \qquad   g \in \{\rho,c_p\}.    
\end{align*}
The model parameters $a_g$ and $b_g$, specified in Table \ref{tab:parameters},  are obtained
through experimental measurements and a corresponding linear least squares fit.

\section{Numerical Solution of Boundary Value Problems} \label{sec: NumericalBVPs}
For numerically solving the fiber boundary value problem the collocation method is embedded into a continuation procedure. The introduced continuation parameter $c\in[0,1]$ regulates the influence of the physical effects in the state equations and boundary conditions, cf., \eqref{eq:num_c}. Proceeding from the solution of a simple auxiliary problem for
$c_0 = 0$, we follow a continuation path $c_j=c_{j-1}+\Delta c_{j-1}$, $j=1,...,n$, and solve a sequence of problems to finally obtain the solution
for our original fiber spinning problem at $c_n=1$.
The core idea of the step size strategy is to compare the effort required to solve the boundary value problem for a whole step $\Delta c$ with the effort required to solve two boundary value problems for the halved step ${\Delta c}/{2}$. The computational effort is measured in evaluations of the right hand side. For details see Algorithm~\ref{algo:adapter}.

As initial step size we choose $\Delta_{c_0} = 0.1$. Furthermore we choose the enlarging factor $\nu_1 = \frac{3}{2}$ and the shrinking factor $\nu_2 = \frac{2}{3}$. 

If a sufficiently good initial guess is available, the direct solution of the BVP with bvp4c is much faster than the solution through Adapter, since in Adapter there have to be solved way more boundary value problems which finally raises computational time. This is also the reason for the outliers with respect to computational time in Table \ref{tab:perfResults}. For parameter $p=(0.3, 11.51)$ the initial guess is not sufficiently good, so bvp4c is not able to determine a solution, and therefore Adapter is used, resulting in higher computational demand.

\IncMargin{1em}
\begin{algorithm}[tb]
\SetKwData{Left}{left}\SetKwData{This}{this}\SetKwData{Up}{up}
\SetKwFunction{Union}{Union}\SetKwFunction{FindCompress}{FindCompress}
\SetKw{Continue}{continue}
\SetKwInOut{Input}{input}\SetKwInOut{Output}{output}
\Input{$f$: parameterized right hand side of ODE\\$g$: parametrized boundary conditions\\ $ y^0$: solution for $c_0=0$\\ $\Delta c_0$: initial step size\\$\mu$: shrinking factor \\ $\nu_1$: shrinking factor\\ $\nu_2$: enlarging factor}
\Output{solution $y^*$ of BVP given by $f$ and $g$}
\BlankLine
$ c_0 \leftarrow 0\,;$ 
$ j \leftarrow 0$ \;
$\Delta c_0 \leftarrow \Delta c_0$\;
\While{$c_j < 1$}{
$\tilde{y}^1 \leftarrow$ solve BVP for $c_j+\Delta c_j $ and initial guess $y^j$\;
$\hat{y}^1  \leftarrow$ solve BVP for $c_j+\Delta c_j/2 $ and initial guess $y^j$\;
$\hat{y}^2  \leftarrow$ solve BVP for $c_j+\Delta c_j $  and initial guess $\hat{y}^1$\;

\If{Newton method diverges}{
$\Delta c_j \leftarrow \frac{\Delta c_j}{\mu} $\;
\Continue
}
$y^* \leftarrow \Tilde{y}^2$\;
\uIf{computation of $\tilde{y}^2$ and $\tilde{y}^3$ need more evaluations of $f$ than computation of $\tilde{y}^1$}{
$\Delta c_{j} \leftarrow \nu_1 \Delta c_j$\;}
\Else{
$ \Delta c_{j} \leftarrow \nu_2 \Delta c_j$\;}
$c_{j+1} \leftarrow \min (1, c_j + \Delta c_j)$\;
$y^{j+1} \leftarrow \tilde{y}^3$ \;
$j \leftarrow j+1$\;
}
\caption{Adapter -- continuation with step size strategy}\label{algo:adapter}
\end{algorithm}\DecMargin{1em}

\section*{Statements and Declarations}
\subsection*{Competing interests}
The authors have no competing interests to declare that are relevant to the content of this article.
\subsection*{Data availability statement}
Data sets generated during the current study are available from the corresponding author on reasonable request.

\bibliographystyle{amsplain}
\bibliography{ref}{}

\providecommand{\bysame}{\leavevmode\hbox to3em{\hrulefill}\thinspace}
\providecommand{\MR}{\relax\ifhmode\unskip\space\fi MR }
\providecommand{\MRhref}[2]{%
  \href{http://www.ams.org/mathscinet-getitem?mr=#1}{#2}
}
\providecommand{\href}[2]{#2}
\begin{thebibliography}{10}

\bibitem{bier2022novel}
Alexander~M. Bier, Walter Arne, and Dirk~W. Schubert, \emph{Novel high-speed
  elongation rheometer}, Macromolecular Materials and Engineering \textbf{307}
  (2022), no.~7, 2100974.

\bibitem{coleman1996interior}
Thomas~F. Coleman and Yuying Li, \emph{An interior trust region approach for
  nonlinear minimization subject to bounds}, SIAM Journal on Optimization
  \textbf{6} (1996), no.~2, 418--445.

\bibitem{hufenus2020melt}
Rudolf Hufenus, Yurong Yan, Martin Dauner, and Takeshi Kikutani,
  \emph{Melt-spun fibers for textile applications}, Materials \textbf{13}
  (2020), no.~19, 4298.

\bibitem{kannengiesserECMI}
Lukas Kannengießer, Nicole Marheineke, and Raimund Wegener, \emph{Generalized
  {N}ewtonian material models in fiber spinning simulations}, European
  Consortium for Mathematics in Industry, Springer, {2024, to appear}.

\bibitem{karian2003handbook}
Harutun Karian, \emph{Handbook of polypropylene and polypropylene composites,
  revised and expanded}, CRC press, 2003.

\bibitem{kierzenka2001bvp}
Jacek Kierzenka and Lawrence~F. Shampine, \emph{A {BVP} solver based on
  residual control and the {M}atlab {PSE}}, ACM Transactions on Mathematical
  Software (TOMS) \textbf{27} (2001), no.~3, 299--316.

\bibitem{kunzelmann2017korrelation}
Peter~Herbert Kunzelmann, \emph{Korrelation der rheologischen {E}igenschaften
  von {M}ischungen aus {P}olypropylenen unterschiedlicher molekularer
  {S}truktur mit deren {V}erarbeitungsverhalten im {S}pinnprozess}, Ph.D.
  thesis, FAU Erlangen-Nürnberg, 2017.

\bibitem{marheineke2011modeling}
Nicole Marheineke and Raimund Wegener, \emph{Modeling and application of a
  stochastic drag for fibers in turbulent flows}, International Journal of
  Multiphase Flow \textbf{37} (2011), no.~2, 136--148.

\bibitem{munstedt2018extensional}
Helmut M{\"u}nstedt, \emph{Extensional rheology and processing of polymeric
  materials}, International Polymer Processing \textbf{33} (2018), no.~5,
  594--618.

\bibitem{qin2019simple}
Yijing Qin and Dirk~W Schubert, \emph{Simple model to predict the effect of
  take-up pressure on fibre diameter of {PET} melt spinning}, Polymer
  \textbf{181} (2019), 121769.

\bibitem{squire1998}
William Squire and George Trapp, \emph{Using complex variables to estimate
  derivatives of real function}, SIAM Review \textbf{40} (1998), no.~1,
  110--112.

\bibitem{sucker1976steady}
Dietrich Sucker and Heinz~P. Brauer, \emph{Steady mass and heat transfer from
  transverse cylinders in steady flow}, W{\"a}rme-und Stoff{\"u}bertragung
  \textbf{9} (1976), 1--12.

\bibitem{trouton1906coefficient}
Frederick~Thomas Trouton, \emph{On the coefficient of viscous traction and its
  relation to that of viscosity}, Proceedings of the Royal Society of London.
  Series A, Containing Papers of a Mathematical and Physical Character
  \textbf{77} (1906), no.~519, 426--440.

\bibitem{wieland2019}
Manuel Wieland, Walter Arne, Robert Fessler, Nicole Marheineke, and Raimund
  Wegener, \emph{An efficient numerical framework for fiber spinning scenarios
  with evaporation effects in airflows}, Journal of Computational Physics
  \textbf{384} (2019), 326--–348.

\end{thebibliography}

\end{document}